\documentclass [11pt, twoside, reqno] {amsart}

\usepackage {geometry}
\usepackage {amsfonts, amssymb, amsmath, amsthm, amscd}
\usepackage {mathrsfs, dsfont, stmaryrd}
\usepackage {stackrel}
\usepackage {color}
\usepackage [colorlinks=true, linkcolor=red, citecolor=red, urlcolor=blue] {hyperref}

\theoremstyle {plain}
\newtheorem {theoreme} {Theorem} [section]
\newtheorem {definition} [theoreme] {Definition}
\newtheorem {proposition} [theoreme] {Proposition}

\newtheorem {lemme} [theoreme] {Lemma}

\newtheorem {remarque} [theoreme] {Remark}

\newtheorem* {propriete*} {Property}
\newtheorem* {proprietes*} {Properties}
\newtheorem* {thanks*} {Acknowledgements}

\numberwithin {equation} {section}

\makeatletter
\def \thebibliographie#1 {\section*{Bibliographie} \list
  {[\arabic{enumi}]}{\settowidth \labelwidth{[#1]} \leftmargin \labelwidth
  \advance \leftmargin \labelsep
  \usecounter{enumi}}
  \def \newblock{\hskip 0.0cm plus 0.0cm minus 0.0cm}
  \sloppy \clubpenalty4000 \widowpenalty4000
  \sfcode`\. = 1000 \relax}

\makeatother

\begin {document}

\title [] {On the lack of compactness on stratified Lie groups}
\author [] {Chieh-Lei Wong}
\address [] {LAMA UMR 8050 \\ Universit\'e Paris-Est Cr\'eteil \\ 61 avenue du G\'en\'eral de Gaulle \\ 94010 Cr\'eteil Cedex}
\email [] {\href {mailto:shell_intheghost@hotmail.com}{shell\_intheghost@hotmail.com}}
\urladdr [] {}
\keywords {\mbox{Lie groups, lack of compactness of the critical Sobolev injection, profiles, wavelets, unconditional bases}}
\date {\today}

\begin {abstract}
\noindent In $\mathbb{R}^d$, the characterization of the \mbox{lack of compactness of the continuous Sobolev injection $ \mathring{H}^s \hookrightarrow L^p $}, with $ \displaystyle{\frac{s}{d} + \frac{1}{p} = \frac{1}{2}} $ and $\displaystyle{0<s<\frac{d}{2}}$, can be rephrased as : a bounded nonzero sequence $f_n$ in $\mathring{H}^s$ admits a subsequence which can be decomposed as a sum of pairwise orthogonal $h$-oscillatory components - known as profiles - and a remainder term which is going strongly to $0$ in $L^p$. The aim of this paper is to generalize this description due to Patrick G\'erard~\cite{Gerard2} to stratified Lie groups. We shall obtain the collection of profiles from a wavelet decomposition by embracing the similar conceptual approach as in~\cite{BCK} or~\cite{Jaffard}.
\end {abstract}
\maketitle

\section {Introduction}
\subsection {Lack of compactness}
\noindent \newline The method of profile decomposition was first  introduced by Ha\"im Br\'ezis and Jean-Michel Coron~\cite{Brezis}, \cite{BC} (see also Michael Struwe~\cite{Struwe}), with roots in the concentrated compactness method of Pierre-Louis Lions~\cite{Lions1}, \cite{Lions2}.
In a paper by Patrick G\'erard~\cite{Gerard2}, the defect of compactness of the Sobolev embedding $ \mathring{H}^s \hookrightarrow L^p $ is described in terms of a sum of rescaled and translated orthogonal profiles, up to a small term in $L^p$. This was generalized to other Sobolev spaces by St\'ephane Jaffard~\cite{Jaffard}, to Besov spaces by Gabriel Koch~\cite{Koch}, and finally to general critical embeddings $ X \hookrightarrow Y $ including a wide range of functional spaces (Lebesgue $L^p$, Sobolev $\mathring{H}^s$, Besov $\mathring{B}^s_{p,q}$, Triebel-Lizorkin $\mathring{F}^s_{p,q}$ only to name a few) by Hajer Bahouri, Albert Cohen and Gabriel Koch~\cite{BCK}. The interested reader can also refer to Kyril Tintarev and Karl-Heinz Fieseler~\cite{TF} for an abstract, functional analytic presentation of the concept in various settings. \\
\newline
\mbox{These profile decomposition techniques have been successfully used for studying nonlinear PDEs, namely :}
\begin {itemize}
\item[\textbullet] the description of bounded energy sequences of solutions of the defocusing semi-linear quintic wave equation, up to remainder terms negligible in energy norm by Hajer Bahouri and Patrick G\'erard~\cite{BGe},
\item[\textbullet] the characterization of the defect of compactness for Strichartz estimates for the Schr\"odinger equation by Sahbi Keraani~\cite{Keraani},
\item[\textbullet] the understanding of features of solutions of nonlinear wave equations with exponential growth by Hajer Bahouri, Mohamed Majdoub and Nader Masmoudi~\cite{BMM},
\item[\textbullet] the sharp estimates of the lifespan of the focusing critical semi-linear wave equation by means of the energy size of the Cauchy data by Carlos E. Kenig and Frank Merle~\cite{KM},
\item[\textbullet] the study of bilinear Strichartz estimates for the wave equation by Terence Tao~\cite{Tao}.
\end {itemize}
For more applications, refer for instance to Isabelle Gallagher and Patrick G\'erard~\cite{GG}, \mbox{Isabelle Gallagher~\cite{Gallagher},} Camille Laurent~\cite{Laurent}, Hajer Bahouri and Isabelle Gallagher~\cite{BGa}, \mbox{Hajer Bahouri,} Jean-Yves Chemin and Isabelle Gallagher~\cite{BCG}, and multiple references therein. \\
\newline
Apart Jamel Benameur~\cite{Benameur} who solved the case of the Heisenberg group (2008), publications were focused on the Euclidean case. This paper provides a positive answer to the natural question of extending the description of the lack of compactness in Sobolev embeddings to stratified Lie groups. Though we will further see a more precise definition, here are some examples of such Lie groups :
\begin {itemize}
\item[\textbullet] Euclidean cases : Abelian fields such as $\mathbb{R}^d$ or $\mathbb{C}^d$, (non Abelian) upper triangular groups,
\item[\textbullet] non-flat cases : Heisenberg groups $\mathbb{H}^d$, Carnot groups, Lie groups of polynomial growth.
\end {itemize}
\noindent We shall also assume that our stratified Lie groups $G$ have a Hausdorff geometrical realization. Any of the above examples are likewise.
\subsection {Stratified Lie groups}
\begin {definition}[FS~\cite{FS}]
A Lie group $(G,\cdot)$ is called stratified if it is connected, simply connected and its Lie algebra $\mathfrak{g}$ decomposes as a direct sum $ \mathfrak{g} = V_1 \oplus \ldots \oplus V_m $ with :
\begin {eqnarray*}
\left\{
\begin {array}{l}
[V_1, V_k] = V_{k+1} \text { if } 1 \leqslant k < m \\
\lbrack V_1, V_m] = \{ 0 \}
\end {array} \right. \; \text{.}
\end {eqnarray*}
\end {definition}
\noindent Thus $\mathfrak{g}$ is a $m$-step, nilpotent and finitely generated, as a Lie algebra, by the vector subspace $V_1$. So as a manifold, $G$ possesses a sub-Riemannian structure. The exponential map is a diffeomorphism from $\mathfrak{g} \longrightarrow G$. When $G$ is identified with $\mathfrak{g}$ via $\exp$, the group law on $G$ is a polynomial map provided by the Campbell-Baker-Hausdorff formula. As a result, the Lie correspondence endows $G$ with a richer manifold structure, and subsequently all the classic notions of differentiable functions, Haar measure, functional spaces, and so on. The left-invariant Haar measure $\mu_G$ on $G$ is induced by the Lebesgue measure on its Lie algebra $\mathfrak{g}$, and we then define the Lebesgue spaces on $G$ as :
\begin {equation*}
L^p(G) = \Bigg\{ \text{Borel functions } f \; \bigg| \; \bigg( \int_G |f|^p d\mu_G \bigg)^{\frac{1}{p}} < +\infty \Bigg\} \; \text{,}
\end {equation*}
with the standard modification when $p=+\infty$. In particular, if left translations are defined by $ \tau_{x'}(x) = x' \cdot x $, the property :
\begin {equation}
\forall f \in L^1(G), \; \forall y \in G, \; \int_G f(y \cdot x) d\mu_G(x) = \int_G f(x) dx \label{left_invariance}
\end {equation}
results from the left-invariance of $\mu_G$. As a homogeneous group, there is a natural action of dilations $(\delta_\alpha)_{\alpha \in \mathbb{R}_+}$ on elements of $G$, given by :
\begin {eqnarray*}
\forall \alpha \in \mathbb{R}_+^*, \; \delta_{\alpha}(x) & = & \alpha \odot x \\
& = & (\underbrace{\alpha x_{1,1}, \ldots, \alpha x_{1, \dim V_1}}_{\begin {subarray}{c} \text{induced by the} \\ \text{canonical action on } V_1 \end {subarray}}, \underbrace{\alpha^2 x_{2,1}, \ldots, \alpha^2 x_{2, \dim V_2}}_{}, \ldots, \underbrace{\alpha^m x_{m,1}, \ldots, \alpha^m x_{m, \dim V_m}}_{}) \; \text{.}
\end {eqnarray*}
In particular, it is immediate to see that :
\begin {eqnarray*}
\delta_1(x) & = & 1_{\mathbb{R}} \odot x = x \; \text{,} \\
\text{and } \big( \delta_{\alpha}(x) \big)^{-1} & = & (\alpha \odot x)^{-1} = \alpha \odot x^{-1} \; \text{.}
\end {eqnarray*}
This family of nonisotropic dilations is a subgroup of $Aut(G)$. Given a dilation $\delta_{\alpha}$, a linear differential operator $X$ is said homogeneous of degree $\ell$ if for any function $f$ on $G$, $ X(f \circ \delta_\alpha) = \alpha^\ell (Xf) \circ \delta_\alpha $. \\
\newline
Since we shall be dealing with a homogeneous group $G$ endowed with a natural family of dilations, we rather define a homogeneous norm $|\cdot|_G$ : $ G \longrightarrow [0,+\infty] $ which is a $\mathscr{C}^{\infty}$ function on $G \backslash \{0\}$ such that :
\begin {equation*}
\left\{
\begin {array}{l}
\forall x \in G, \; |x^{-1}|_G = |x|_G \\
\forall \alpha > 0, \; | \alpha \odot x |_G = \alpha |x|_G \\
|x|_G = 0 \text{ iff } x=0
\end {array} \right. \; \text{.}
\end {equation*}
Let us denote by $ \displaystyle{Q = \sum_{k=1}^m k \dim_{\mathbb{R}} V_k} $ the homogeneous dimension of $G$, and introduce the $L^1$-normalized dilation $\underline{\delta_t^1} f$ of a function $f$ by :
\begin {equation}
\forall f, \; \forall x \in G, \; \forall t>0, \; \underline{\delta_t^1} f(x) = t^Q f(t \odot x) \; \text{.}
\end {equation}
Note that, by definition, $\underline{\delta_t^1}$ preserves the norm $L^1$, that is $ \| \underline{\delta_t^1} f \|_{L^1(G)} = \| f \|_{L^1(G)} $.
\begin {remarque}
Depending on the context i.e. the functional spaces we are working with, we might use different normalizations for a given function. The superscript in $\underline{\delta_t^{\bullet}}$ is helpful to remind which normalization is chosen.
\end {remarque}
\noindent Elements of $\mathfrak{g}$ can be identified with differential operators of length 1 of $G$ which are invariant under left translations. Note that a vector field $X$ : $G \longrightarrow TG$ is said to be left-invariant when the following diagram commutes for all $h \in G$ :
\begin {equation*}
\begin {CD}
G @>\tau_h>> G \\
@VXVV @VVXV\\
TG @>>d\tau_h> TG
\end {CD}
\end {equation*}
where $\tau_h$ is the left translation on $G$ defined by $\tau_h(x) = h \cdot x$. Then it follows that for all $h$ in $G$, $\; X \circ \tau_h = d\tau_h \circ X$. This infinitesimal characterization is equivalent to say that for any smooth function $f$, one has $X(f \circ \tau_h) = (Xf) \circ \tau_h$. \\
\newline
Let $n \in \mathbb{N}$. Let $I$ and $k$ be multi-indexes in $\mathbb{N}^n$. For a differential operator $ X^I = X_{k_1}^{I_1} X_{k_2}^{I_2} \ldots X_{k_n}^{I_n} $ where the $X_{k_i}$'s are taken in $\mathfrak{g}$, we define two distinctive notions :
\begin {itemize}
\item[\textbullet] its isotropic length (or order) : $|I| = I_1 + I_2 + \ldots + I_n$,
\item[\textbullet] its homogeneous degree : $ \displaystyle{\deg X^I = \sum_{i=1}^n I_i \deg X_{k_i}} $, where $\deg X_{k_i} = j$ if $X_{k_i} \in V_j$.
\end {itemize}
\noindent \newline Under the identification $\mathfrak{g} \simeq G$, polynomials on $G$ are polynomials on $\mathfrak{g}$. Let us denote by $\mathcal{P}_k$ the vector space of homogeneous polynomials of degree $k$, and set $ \mathcal{P} = \varinjlim \mathcal{P}_k $. One can also define the Schwartz space on $G$ by $ \mathcal{S}(G) = \mathcal{S}(\mathfrak{g}) $. Let $ \mathcal{S}'(G) $ be the space of tempered distributions on $G$ and $ \mathcal{S}'(G) / \mathcal{P} $ the space of tempered distributions modulo polynomials on $G$. Duality between the two spaces is achieved by the sesquilinear product $\langle \cdot,\cdot \rangle$ : $ \mathcal{S}'(G) \times \mathcal{S}(G) \longrightarrow \mathbb{C} $ defined by $ \langle f,g \rangle = \displaystyle{\int_G f \bar{g} \; d\mu_G} $. \\
\newline
It is common to use the spectral calculus of a suitable sub-Laplacian :
\begin {equation}
\Delta_G = \sum_{X_j \in V_1} X_j^2 \; \text{,}
\end {equation}
induced by the aforementioned sub-Riemannian structure of $G$ in the intent to define a Littlewood-Paley decomposition for functions and tempered distributions on $G$. Restricted to $\mathscr{C}_c^\infty(G)$ - that is the space of smooth functions with compact support defined on $G$, the sub-Laplacian $\Delta_G$ is a linear differential operator, homogeneous of degree 2 and formally self-adjoint :
\begin {equation*}
\forall f,g \in \mathscr{C}_c^\infty(G), \; \langle \Delta_G f,g \rangle = \langle f,\Delta_G g \rangle \; \text{.}
\end {equation*}
Its closure has a domain $ \mathscr{D} = \Big\{ u \in L^2(G) \; \big| \; \Delta_G u \in L^2(G) \Big\} $ where $\Delta_G u$ is taken in the sense of distributions. It follows that its closure is also self-adjoint and it is actually the unique self-adjoint extension of $\Delta_{|\mathscr{C}_c^\infty(G)}$. We shall continue to denote this extension by $\Delta_G$. \\
Moreover, the homogeneous Sobolev spaces $ \mathring{H}^s(G) $ are defined as :
\begin {equation}
\mathring{H}^s(G) = \bigg\{ f \in \mathcal{S}'(G) / \mathcal{P} \; \Big| \; (-\Delta_G)^{\frac{s}{2}} f \in L^2(G) \bigg\} \; \text{.}
\end {equation}
\subsection {Characterization of the lack of compactness of the Sobolev injection $\mathring{H}^s(G) \hookrightarrow L^p(G)$}
\noindent \newline Set $ \displaystyle{\frac{s}{Q} + \frac{1}{p} = \frac{1}{2}} $ where $\displaystyle{0<s<\frac{Q}{2}} \; \cdot $ In the framework of the Sobolev injection (see e.g.~\cite{BFG}, or~\cite{BGX} for the  specific case of the Heisenberg group $\mathbb{H}^d$) :
\begin {equation} \label{critical_sobolev_injection}
\mathring{H}^s(G) \hookrightarrow L^p(G) \; \text{,}
\end {equation}
both spaces are homogeneous spaces with the same scaling properties since :
\begin {eqnarray*}
& h^{\frac{Q}{p}} \| f \circ \delta_h \|_{\mathring{H}^s(G)} = \|f\|_{\mathring{H}^s(G)} \; \text{,} \\
\text{and} & h^{\frac{Q}{p}} \| f \circ \delta_h \|_{L^p(G)} = \|f\|_{L^p(G)} \; \text{.}
\end {eqnarray*}
For any function $u \in \mathring{H}^s(G)$, if we define the operators :
\begin {eqnarray*}
& & \underline{\tau_\kappa} u = u \circ \tau_{\kappa^{-1}} \; \text{,} \\
& \text{and} & \underline{\delta_h^p} u = h^{\frac{Q}{p}} u \circ \delta_h \; \text{,}
\end {eqnarray*}
both $L^p$ and $\mathring{H}^s$ norms are preserved under translations $u \longmapsto \underline{\tau_\kappa} u$ and dilations $ u \longmapsto \underline{\delta_{h^{-1}}^p} u $. If $u$ is a nonzero element of $\mathring{H}^s(G)$, for any sequence of points $ (\kappa_n)_{n \in \mathbb{N}} $ going to infinity in $G$ (i.e $ | \kappa_n |_G \stackrel[n \to +\infty]{}{\longrightarrow} +\infty $) and any sequence of positive real numbers $ (h_n)_{n \in \mathbb{N}} $ converging to $0$ or $+\infty$, the two sequences $ \Big( \underline{\tau_{\kappa_n}} u \Big)_{n \in \mathbb{N}}$ and $ \displaystyle{ \Big( \underline{\delta_{h_n^{-1}}^p} u \Big)_{n \in \mathbb{N}}} $ converge weakly to $0$ in $\mathring{H}^s(G)$, henceforth they are not relatively compact in $L^p(G)$. In this paper, we shall prove that these invariances under $\underline{\tau_{\kappa_n}}$ and $\underline{\delta_{h_n^{-1}}^p}$ are single responsible for the lack of compactness of the continuous Sobolev injection $ \mathring{H}^s(G) \hookrightarrow L^p(G) $. \\
\newline
Before stating the main result, let us introduce the notions of scales and concentration cores.
\begin {definition} \label{scales_cores}
We call a scale any sequence $\underline{h}=(h_n)_{n \in \mathbb{N}}$ of positive real numbers, and a concentration core any sequence $\underline{\kappa}=(\kappa_n)_{n \in \mathbb{N}}$ of points in $G$. Pairs $(\underline{h},\underline{\kappa})$ and $(\underline{\tilde h},\underline{\tilde \kappa})$ are said orthogonal if :
\begin {eqnarray}
\log \bigg| \frac{h_{n}}{\tilde h_{n}} \bigg| \stackrel[n \to +\infty]{}{\longrightarrow} \pm \infty & & \text{for the scales,} \label{echelles_etrangeres} \\
\text{or } \bigg( h_{n} = \tilde h_{n} \text{ and } \frac{1}{h_{n}} | \kappa_{n}^{-1} \cdot \tilde \kappa_{n} |_G \stackrel[n \to +\infty]{}{\longrightarrow} +\infty \bigg) & & \text{for the concentration cores.}
\end {eqnarray}
\end {definition}
\begin {theoreme} \label{theoreme_profils}
Consider the continuous embedding $ \mathring{H}^s(G) \hookrightarrow L^p(G) $ with $ \displaystyle{\frac{s}{Q} + \frac{1}{p} = \frac{1}{2}} $ and $ \displaystyle{0 < s < \frac{Q}{2}} \; \cdot $ Let $(u_n)_{n>0}$ be a sequence of bounded functions in $ \mathring{H}^s(G) $. Then, up to the possible extraction of a subsequence, there exist a family of functions $(\phi^{\ell})_{\ell \in \mathbb{N}^*}$ in $ \mathring{H}^s(G) $ - the so-called profiles - as well as families of scales $(\underline{h}^{\ell})=(h^\ell _{n})$ and concentration cores $(\underline{\kappa}^{\ell})=(\kappa^\ell _{n})$ such that :
\begin {itemize}
\item[(i)] the pairs $(\underline{h}^{\ell},\underline{\kappa}^{\ell})$ are pairwise orthogonal in the sense of Definition \ref{scales_cores},
\item[(ii)] for any $L \geqslant 1$, we have :
\begin {eqnarray} \label{profils}
& & u_n(x) = \underbrace{\sum_{\ell=1}^{L} (h^\ell _{n})^{s-\frac{Q}{p}} \phi^{\ell} \bigg( \frac{1}{h^\ell _{n}} \odot ((\kappa^\ell _{n})^{-1} \cdot x) \bigg)}_{\begin {subarray}{c} \text{superposition of } \\ h_{\ell,n} \text{-oscillatory components} \end {subarray}} + r_{n,L}(x) \; \text{,} \\
& \text{with} & \lim_{L \to +\infty} \varlimsup_{n \to +\infty} \| r_{n,L} \|_{L^p(G)} = 0 \; \text{.} \nonumber
\end {eqnarray}
\end {itemize}
The profile decomposition (\ref{profils}) is asymptotically orthogonal (or almost orthogonal) in the sense that :
\begin {equation}
\| u_n \|_{\mathring{H}^s(G)}^2 = \sum_{\ell=1}^L \| \phi^{\ell} \|_{\mathring{H}^s(G)}^2 + \| r_{n,L} \|_{\mathring{H}^s(G)}^2 + o(1) \text{ as } n \to +\infty \; \text{.} \label{almost_orthogonality}
\end {equation}
\end {theoreme}
\noindent This profile decomposition recovers various versions of the concentration-compactness principle, as the ones in~\cite{Gerard2}, \cite{Lions1} or \cite{Lions2}.
\begin {remarque}
For some sequences $(u_n)_{n>0}$, it may happen that the decomposition (\ref{profils}) contains only a finite number of profiles. In particular, the sequence $(u_n)_{n>0}$ is compact in $L^p(G)$ iff $\forall \ell \geqslant 1, \; \phi^{\ell} = 0 $.
\end {remarque}
\subsection {Layout of this paper}
\noindent \newline We shall prove Theorem~\ref{theoreme_profils} by transposing to stratified Lie groups the method developed in~\cite{BCK} which is mainly based on wavelet decompositions. The authors considered critical Sobolev embeddings $ X \hookrightarrow Y $ for generic functional spaces $X$ and $Y$ with same scaling properties and endowed with unconditional wavelet basis. They also emphasized on two key properties : the first one is tied to nonlinear approximation, and the other one is similar to Fatou's lemma. \\
To achieve our goal, we firstly exhibit a wavelet basis by using the spectral calculus of the sub-Laplacian. The mother wavelet $\psi$ will feature some nice properties (in $\mathcal{S}(G)$, with infinite vanishing moments in the space variable, and compactly supported in its conjugate variable). Translated and dilated copies of $\psi$ provide an unconditional wavelet basis for both homogeneous Besov spaces $ \mathring{B}_{p,q}^s(G) $ with $1 \leqslant p,q < +\infty$, and $ L^p(G) $ with $1<p<+\infty$. Consequently, the rest of this paper is divided into three parts :
\begin {itemize}
\item[\textbullet] Section~\ref{LP-admissible_wavelets} deals with Littlewood-Paley (abbreviated as LP for short) -admissible functions $\psi$ in order to characterize the homogeneous Besov spaces $\mathring{B}^s_{p,q}(G)$,
\item[\textbullet] in Section~\ref{besov_spaces}, we will clarify the unconditional convergence of the wavelet expansion in both $L^p(G)$ and $\mathring{B}^s_{p,q}(G)$ via discrete sampling techniques and Banach frames,
\item[\textbullet] in Section~\ref{profils_stratified_groups}, we will perform the algorithm for the profiles' extraction on $G$, and eventually prove Theorem~\ref{theoreme_profils}.
\end {itemize}
\noindent Note that in Sections~\ref{LP-admissible_wavelets} and~\ref{besov_spaces}, we shall make straightforward use of several results from FS~\cite{FS} and FM~\cite{Fuhr}. We intentionally omit some technical proofs which can be easily found in these references.

\begin {thanks*}
\noindent The author would like to thank Clotilde Fermanian-Kammerer and Hajer Bahouri for numerous fruitful discussions.
\end {thanks*}

\section {LP-admissible functions and Besov spaces} \label{LP-admissible_wavelets}
\noindent The purpose of this section is to define the homogeneous Besov spaces $ \mathring{B}^s_{p,q}(G) $ via a Littlewood-Paley decomposition. In the Euclidean framework, a LP-admissible function $\psi$ is constructed by defining a dyadic partition of unity in the Fourier side, and then apply inverse Fourier $\mathscr{F}^{-1}$. The standard method of transposing this construction to stratified Lie groups is to replace the Fourier transform by the spectral decomposition of the sub-Laplacian $\Delta_G$ (see further Lemma~\ref{lemme_fondamental}), and it is essentially based on convolution. Then this construction proves to be independent of the choice of the basis of the stratum $V_1$. The proper framework for such Littlewood-Paley decompositions is actually the previously introduced space of tempered distributions modulo polynomials $ \mathcal{S}'(G) / \mathcal{P} $. The convergence of these decompositions is obtained via duality with the subspace $ \mathcal{V}(G) \subset \mathcal{S}(G) $ of functions with an infinite number of vanishing moments. \\
\newline
Let us recall that the noncommutative convolution product of two functions $f$ and $g$ on $G$ is defined by :
\begin {equation*}
f*g(x) = \int_G f(x \cdot y^{-1}) g(y) d\mu_G(y) = \int_G f(y) g(y^{-1} \cdot x) d\mu_G(y) \; \text{,}
\end {equation*}
and if $X$ is a left-invariant vector field on $G$, then :
\begin {equation*}
X(f*g)(x) = f*Xg(x) = \int_G f(y) (Xg)(y^{-1} \cdot x) d\mu_G(y) \; \text{.}
\end {equation*}
In what follows, for a function $f$ on $G$, we define $ \tilde{f}(x) = f(x^{-1}) $ and set :
\begin {equation}
\forall x \in G, \; f^*(x) = \overline{f(x^{-1})} \Longleftrightarrow f^* = \bar{\tilde{f}} \; \text{.}
\end {equation}
Then one can easily verify that for $f$ in $ L^2(G) \cap L^1(G) $, the adjoint of the convolution operator $ g \longmapsto g*f $ is $ g \longmapsto g*f^* $.
\subsection {The subspace $\mathcal{V}(G)$ of functions with an infinite number of vanishing moments}
\begin {definition}
Let $ n \in \mathbb{N} $, a function $f$ : $ G \longrightarrow \mathbb{C} $ is said of polynomial decay of order $n$ if there exists a constant $c>0$ such that :
\begin {equation*}
\forall x \in G, \; |f(x)| \leqslant \frac{c}{(1+|x|_G)^n} \; \cdot
\end {equation*}
\end {definition}
\begin {definition}
Let $ r \in \mathbb{N} $, a function $f$ has vanishing moments of order $r$ if :
\begin {equation*}
\forall p \in \mathcal{P}_{r-1}, \; \int_G f(x) p(x) dx = 0 \; \text{.}
\end {equation*}
\end {definition}
\noindent Under the identification of $G$ with $\mathfrak{g}$, the inverse map $x \longmapsto x^{-1}$ is identified with the additive inverse map $x \longmapsto -x$ on $\mathfrak{g}$. It follows that $ \forall p \in \mathcal{P}_{r-1}, \; \tilde{p} \in \mathcal{P}_{r-1} $, where $\tilde{p}$ is defined as above. If $f$ has vanishing moments of order $r$ then :
\begin {equation*}
\forall p \in \mathcal{P}_{r-1}, \; \int_G \tilde{f}(x) p(x) dx = \int_G f(x) \tilde{p}(x) dx = 0 \; \text{,}
\end {equation*}
which shows that $\tilde{f}$ has also vanishing moments of order $r$. This property of having vanishing moments is central in the wavelet analysis due to the following principle : in a convolution product such as $ g * \underline{\delta_t^1} f $, vanishing moments of one factor coupled with smoothness and regularity of the other factor lead to a rapidly decay at infinity of $ g * \underline{\delta_t^1} f $. Thereafter, we shall denote by $\mathcal{V}(G)$ the subspace in the Schwartz class $\mathcal{S}(G)$ of functions with an infinite number of vanishing moments. We sum up below some useful topological properties of $\mathcal{S}(G)$ :
\begin {proprietes*}[FS~\cite{FS} - \S 1 D. The Schwartz class]
\noindent
\begin {itemize}
\item[\textbullet] $\mathcal{S}(G)$ is a Fr\'echet space whose topology is conveniently defined by any norm of the family of norms :
\begin {equation*}
\forall n \in \mathbb{N}, \; \| f \|_n = \sup_{\substack{x \in G \\ |I| \leqslant n}} (1+|x|_G)^{(n+1)(Q+1)} \big| X^I f (x) \big| \; \text{.}
\end {equation*}
Convolution is continuous from $\mathcal{S}(G) \times \mathcal{S}(G)$ to $\mathcal{S}(G)$. More precisely, for every $n \in \mathbb{N}$, there exists $c_n > 0$ such that :
\begin {equation*}
\forall \varphi, \psi \in \mathcal{S}, \; \| \varphi * \psi \|_n \leqslant c_n \| \varphi \|_n \| \psi \|_{n+1} \; \text{.}
\end {equation*}
\item[\textbullet] $\mathcal{V}(G)$ is a closed subspace (in particular, it is complete) of $\mathcal{S}(G)$, with $ \mathcal{S}(G) * \mathcal{V}(G) \subset \mathcal{V}(G) $. \\
If $\mathcal{V}'(G)$ denotes the topological dual of $\mathcal{V}(G)$, then $\mathcal{V}'(G)$ can be canonically identified with $ \mathcal{S}'(G) / \mathcal{P} $.
\item[\textbullet] For all $ \psi \in \mathcal{S}(G) $, the map $ \mathcal{S}'(G) / \mathcal{P} \longrightarrow \mathcal{S}'(G) / \mathcal{P} $ defined by $ u \longmapsto u * \psi $ is a well-defined operator and it is continuous on $\mathcal{S}'(G) / \mathcal{P}$. If $ \psi \in \mathcal{V}(G) $, the associated convolution operator is a well-defined and continuous operator from $ \mathcal{S}'(G) / \mathcal{P} $ into $ \mathcal{S}'(G) $.
\end {itemize}
\end {proprietes*}
\subsection {The Calder\'on's reproducing formula}
\noindent \newline From here, we shall use the classic notation $\hat{f} = \mathscr{F}f$ for the Fourier transform on $\mathbb{R}$. Suppose that $\Delta_G$ has a spectral resolution $ \displaystyle{\Delta_G = \int_0^{+\infty} \lambda dP_{\lambda}} $, where $dP_{\lambda}$ is the spectral projection measure i.e. a measure with values in the spectrum of $\Delta_G$. For any bounded Borel function $\hat{f}$ on $\mathbb{R}_+$, the operator $ \displaystyle{\hat{f}(\Delta_G) = \int_0^{+\infty} \hat{f}(\lambda) dP_{\lambda}} $ which should be understood as :
\begin {equation*}
\forall \phi,\eta \in L^2(G), \; \langle \hat{f}(\Delta_G) \phi, \eta \rangle = \int_0^{+\infty} \hat{f}(\lambda) dP_{\lambda}(\phi,\eta)
\end {equation*}
$\big($where $dP_{\lambda}(\phi,\eta)$ is the unique Borel measure associated to the pair $(\phi,\eta) \big)$ is a bounded integral operator in $\mathscr{L}^2$ with a convolution distribution-kernel $K_f$ in $L^2(G)$ satisfying $ \forall \eta \in \mathcal{S}(G), \; \hat{f}(\Delta_G) \eta = \eta * K_f $. In the following, we shall denote this kernel $K_f$ in some abusive way by $f$, that is :
\begin {equation} \label{definition_kernel_fourier}
\forall \eta \in \mathcal{S}(G), \; \hat{f}(\Delta_G) \eta = \eta * f \; \text{.}
\end {equation}
An important property due to Andrzej Hulanicki~\cite{Hulanicki} is that, for smooth and rapidly decaying functions $ \hat{f} \in \mathcal{S}(\mathbb{R}_+) $, the kernel associated to $ \hat{f}(\Delta_G) $ is a function in the Schwartz class, namely $f \in \mathcal{S}(G)$. \\
\newline
At this stage, we have all the ingredients for a Littlewood-Paley type decomposition :
\begin {equation} \label{calderon_reproducing_formula}
f = \sum_{j \in \mathbb{Z}} f * \psi_j^* * \psi_j \; \text{,}
\end {equation}
where $\forall j \in \mathbb{Z}, \; \psi_j = \underline{\delta_{2^j}^1} \psi$ is a dilated of $ \psi \in \mathcal{S}(G) $ (with the usual convention $ \psi_0 = \psi $).
\begin {definition}
A function $ \psi \in \mathcal{S}(G) $ is said LP-admissible if :
\begin {equation*}
\forall g \in \mathcal{V}(G), \; g = \lim_{n \to +\infty} \sum_{|j| \leqslant n} g * \psi_j^* * \psi_j
\end {equation*}
holds, with convergence in $\mathcal{S}(G)$. Duality induces the convergence of the decomposition on $ \mathcal{S}'(G) / \mathcal{P} $ :
\begin {equation*}
\forall u \in \mathcal{S}'(G) / \mathcal{P}, \; u = \lim_{n \to +\infty} \sum_{|j| \leqslant n} u * \psi_j^* * \psi_j \; \text{.}
\end {equation*}
\end {definition}
\noindent Before stating Lemma~\ref{lemme_fondamental} which is the cornerstone for the construction of LP-admissible functions, we need the next useful preliminary result.
\begin {lemme}[FS~\cite{FS} - Proposition 1.49] \label{lemme_limite_at_infty}
Let $\vartheta \in \mathcal{S}(G)$ and define $ \vartheta_t = \underline{\delta_t^1} \vartheta $, then :
\begin {eqnarray*}
& & \forall \psi \in \mathcal{S}(G), \; \psi * \vartheta_{\frac{1}{t}} \longrightarrow c_{\vartheta} \psi \text{ in } \mathcal{S}(G) \; \text{,} \\
& \text{and} & \forall f \in \mathcal{S}'(G), \; f * \vartheta_{\frac{1}{t}} \longrightarrow c_{\vartheta} f \text{ in } \mathcal{S}'(G) \; \text{,}
\end {eqnarray*}
as $t \to 0$, where $\displaystyle{c_{\vartheta} = \int_G \vartheta(x) dx}$.
\end {lemme}
\begin {lemme} \label{lemme_fondamental}
Let $\hat{\phi}$ be a function in $\mathscr{C}^\infty(\mathbb{R})$ with support in $[0,4]$ such that $ \displaystyle {\left\{
\begin {array}{l}
0 \leqslant \hat{\phi} \leqslant 1 \\
\hat{\phi} \equiv 1 \text{ on } \displaystyle{\bigg[0,\frac{1}{4}\bigg]}
\end {array} \right.} $. \\
Set $ \hat{\psi}(\xi) = \sqrt{\hat{\phi}(2^{-2}\xi) - \hat{\phi}(\xi)} $. Thus $ \hat{\psi} \in \mathscr{C}_c^\infty(\mathbb{R}_+) $ with support in the interval $ \displaystyle{\bigg[ \frac{1}{4},4 \bigg]} $ and we get a dyadic partition of unity with $\hat{\psi}$ : $ \displaystyle{\sum_{j \in \mathbb{Z}} \big| \hat{\psi}(2^{-2j}\xi) \big|^2 = 1} $ a.e. \\
Let $\Delta_G$ be the sub-Laplacian. Let $\psi$ be the convolution distribution-kernel as defined in (\ref{definition_kernel_fourier}) which is associated to the bounded left-invariant operator $\hat{\psi}(\Delta_G)$, then $\psi$ is LP-admissible and belongs to $\mathcal{V}(G)$.
\end {lemme}
\begin {proof}
The spectral theorem applied to the dyadic partition of unity with $\hat{\psi}$ gives :
\begin {equation*}
\sum_{j \in \mathbb{Z}} \bigg[ \hat{\psi}(2^{-2j} \Delta_G) \bigg]^* \circ \bigg[ \hat{\psi}(2^{-2j} \Delta_G) \bigg] = \mathds{1} \; \text{.}
\end {equation*}
Let $g \in \mathcal{V}(G)$. First observe that due to the quadratic homogeneity of $\Delta_G$, the convolution kernel associated to $\hat{\psi}(2^{-2j} \Delta_G)$ coincides with $\psi_j = \underline{\delta_{2^j}^1} \psi$. The decomposition :
\begin {equation} \label{equation_spectrale}
g = \sum_{j \in \mathbb{Z}} \bigg[ \hat{\psi}(2^{-2j} \Delta_G) \bigg]^* \circ \bigg[ \hat{\psi}(2^{-2j} \Delta_G) \bigg] g = \sum_{j \in \mathbb{Z}} g * \psi_j^* * \psi_j
\end {equation}
holds in norm $L^2$. Note that, since $\hat{\psi}$ is a real-valued function, one has actually $ \psi^* = \psi $. \\
For any integer $m \in \mathbb{N}$, one has :
\begin {equation*}
\sum_{|j| \leqslant m} g * \psi_j^* * \psi_j = g * \underline{\delta_{2^{m+1}}^1} \phi - g * \underline{\delta_{2^{-m}}^1} \phi \; \text{,}
\end {equation*}
where $\phi \in \mathcal{S}(G)$ is the convolution kernel of $\hat{\phi}(\Delta_G)$. Since $\phi$ is in the Schwartz class, it follows by Lemma~\ref{lemme_limite_at_infty} that $ \displaystyle{g * \underline{\delta_{2^{m+1}}^1} \phi \stackrel[m \to +\infty]{}{\longrightarrow} c_{\phi} g} $ in $\mathcal{S}(G)$ for some constant $c_{\phi}$. Hence $ \displaystyle{\sum_{|j| \leqslant m} g * \psi_j^* * \psi_j \stackrel[m \to +\infty]{}{\longrightarrow} c_{\phi} g} $ in $\mathcal{S}(G)$, and the identity (\ref{equation_spectrale}) in $L^2(G)$ gives $c_{\phi}=1$.
\end {proof}
\noindent The Calder\'on's decomposition $ \displaystyle{g = \sum_{j \in \mathbb{Z}} g * \psi_j^* * \psi_j} $ converges strongly and unconditionally in norm $L^2$ because of the unconditional convergence of the sum $ \displaystyle{\sum_{j \in \mathbb{Z}} \overline{\hat{\psi}(2^{-2j} \xi)} \hat{\psi}(2^{-2j} \xi)} $.
\begin {proprietes*}
\noindent
\begin {itemize}
\item[\textbullet] $ \psi \in \mathcal{S}(G) $~\footnote{Andrzej Hulanicki : A functional calculus for Rockland operators on nilpotent Lie groups (1984)} and has an infinite number of vanishing moments~\footnote{Daryl Geller, Azita Mayeli : Continuous wavelets and frames on stratified Lie groups (2006)}, hence $ \psi \in \mathcal{V}(G) $.
\item[\textbullet] Any function LP-admissible built according to Lemma~\ref{lemme_fondamental} satisfies the relation :
\begin {equation} \label{relation_rangs}
\forall j,\ell \in \mathbb{Z}, \; |j-\ell|>1 \; \Longrightarrow \; \psi_j^* * \psi_\ell = 0 \; \text{,}
\end {equation}
resulting from :
\begin {equation*}
\bigg[ \hat{\psi}(2^{-2j} \Delta_G) \bigg]^* \circ \bigg[ \hat{\psi}(2^{-2\ell} \Delta_G) \bigg] = 0 \; \text{.}
\end {equation*}
\end {itemize}
\end {proprietes*}
\subsection {Homogeneous Besov spaces}
\begin {definition} \label{definition_besov}
Let $ \psi \in \mathcal{V}(G) $ be LP-admissible. Let $ 1 \leqslant p,q \leqslant +\infty $ and $ s \in \mathbb{R} $. The homogeneous Besov space associated to $\psi$ is defined by :
\begin {equation}
\mathring{B}^{s,\psi}_{p,q}(G) = \bigg\{ u \in \mathcal{S}'(G) / \mathcal{P} \; \Big| \; \big\{ 2^{js}\| u * \psi_j^* \|_{L^p(G)} \big\}_{j \in \mathbb{Z}} \in \ell^q(\mathbb{Z}) \bigg\} \; \text{,}
\end {equation}
with the associated norm :
\begin {equation}
\| u \|_{\mathring{B}^{s,\psi}_{p,q}(G)} = \Big\| \big\{ 2^{js}\| u * \psi_j^* \|_{L^p(G)} \big\}_{j \in \mathbb{Z}} \Big\|_{\ell^q(\mathbb{Z})} \; \text{.}
\end {equation}
\end {definition}
\noindent The definition of homogeneous Besov spaces requires taking $L^p$ norms of elements of $ \mathcal{S}'(G) / \mathcal{P} $. We use the canonical embedding $ L^p(G) \hookrightarrow \mathcal{S}'(G) $. For $p<+\infty$, by using that $ \mathcal{P} \cap L^p(G) = \{0\} $, one has the embedding $ L^p(G) \hookrightarrow \mathcal{S}'(G) / \mathcal{P} $. Given $ u \in \mathcal{S}'(G) / \mathcal{P} $, we define :
\begin {itemize}
\item[\textbullet] $ \| u \|_{L^p(G)} = \| u+q \|_{L^p(G)} $ when $ u+q \in L^p(G) $, for a suitable $ q \in \mathcal{P} $ ; note that the decomposition $u+q$ is unique since $ \mathcal{P} \cap L^p(G) = \{0\} $,
\item[\textbullet] $ \| u \|_{L^p(G)} = +\infty $ otherwise.
\end {itemize}
By contrast, the norm $ \| \cdot \|_{L^\infty(G)} $ can only be defined on $ \mathcal{S}'(G) $ by assigning the value $+\infty$ to any \mbox{$ u \in \mathcal{S}'(G) \backslash L^\infty(G) $}. Note that the Hausdorff-Young's inequality $ \| u*f \|_{L^p(G)} \leqslant \| u \|_{L^p(G)} \| f \|_{L^1(G)} $ remains valid respectively :
\begin {eqnarray*}
\text{for } p<+\infty & \text{:} & \forall f \in \mathcal{S}(G), \; \forall u \in \mathcal{S}'(G) / \mathcal{P} \; \text{,} \\
\text{and for } p=+\infty & \text{:} & \forall f \in \mathcal{S}(G), \; \forall u \in \mathcal{S}'(G) \; \text{.}
\end {eqnarray*}
For $p<+\infty$, if $ u+q \in L^p(G) $ then $ (u+q) * \psi = u * \psi + q * \psi \in L^p(G) $. \\
\newline The combination of Lemma~\ref{lemme_fondamental} and Definition~\ref{definition_besov} shows that we shall recover the usual notion of homogeneous Besov spaces built via the spectral calculus of sub-Laplacians. The definition turns out to be independent of the choice of $\psi$.
\begin {theoreme}[FM~\cite{Fuhr} - Theorem 3.11] \label{psi_independent}
Let $ \psi_1, \psi_2 \in \mathcal{V}(G) $ be LP-admissible. Let $ 1 \leqslant p,q \leqslant +\infty $ and $ s \in \mathbb{R} $. Then $ \mathring{B}^{s,\psi_1}_{p,q}(G) = \mathring{B}^{s,\psi_2}_{p,q}(G) $ with equivalence of norms.
\end {theoreme}
\noindent \mbox{We will accordingly omit the $\psi$ superscript and simply write $\mathring{B}^s_{p,q}(G)$ for any choice of LP-admissible $ \psi \in \mathcal{V}(G) $.}
\begin {proprietes*}
\noindent
\begin {itemize}
\item[\textbullet] $\mathring{B}^s_{p,q}(G)$ is a Banach space.
\item[\textbullet] For $ 1 \leqslant p,q \leqslant +\infty $ and $ s \in \mathbb{R} $, one has the following continuous embeddings :
\begin {eqnarray*}
& & \mathcal{V}(G) \hookrightarrow \mathring{B}^s_{p,q}(G) \hookrightarrow \mathcal{S}'(G) / \mathcal{P} \; \text{,} \\
\text{and} & & \mathcal{V}(G) \hookrightarrow \big( \mathring{B}^s_{p,q}(G) \big)^{'} = \mathring{B}^{-s}_{\bar{p},\bar{q}}(G) \; \text{.}
\end {eqnarray*}
For $ p,q < +\infty $, $\mathcal{V}(G)$ is dense in $\mathring{B}^s_{p,q}(G)$.
\item[\textbullet] The space $ \displaystyle{\mathcal{S}(G) * \sum_{|j| \leqslant n} \psi_j^* * \psi_j} $ is dense in $\mathcal{V}(G)$, as well as in $\mathring{B}^s_{p,q}(G)$ if $ p,q < +\infty $. And in this space, the decomposition $ \displaystyle{g = \sum_{j \in \mathbb{Z}} g * \psi_j^* * \psi_j} $ has a finite number of nonzero terms.
\end {itemize}
\end {proprietes*}
\noindent We now extend this Littlewood-Paley decomposition to other functional spaces.
\begin {proposition}[FM~\cite{Fuhr} - Proposition 3.14]
Let $ 1 \leqslant p,q < +\infty $ and $ \psi \in \mathcal{V}(G) $ be LP-admissible. Then for all $g$ in $\mathring{B}^s_{p,q}(G)$ :
\begin {equation*}
g = \lim_{n \to +\infty} \sum_{|j| \leqslant n} g * \psi_j^* * \psi_j
\end {equation*}
holds in $\mathring{B}^s_{p,q}(G)$.
\end {proposition}
\begin {proof}
Consider the operators $\Sigma_n$ : $ \mathring{B}^s_{p,q}(G) \longrightarrow \mathring{B}^s_{p,q}(G) $ defined by :
\begin {equation*}
\Sigma_n g = \sum_{|j| \leqslant n} g * \psi_j^* * \psi_j \; \text{.}
\end {equation*}
This family of operators $(\Sigma_n)_{n \in \mathbb{N}}$ is bounded in norm. The $\Sigma_n$'s converge strongly to the identity operator on a dense subspace of $\mathring{B}^s_{p,q}(G)$. But by boundedness of the family, this implies the strong convergence everywhere.
\end {proof}
\begin {proposition}[FM~\cite{Fuhr} - Proposition 3.15]
Let $ 1 < p < +\infty $ and $ \psi \in \mathcal{V}(G) $ be LP-admissible. Then for all $g$ in $L^p(G)$ :
\begin {equation*}
g = \lim_{n \to +\infty} \sum_{|j| \leqslant n} g * \psi_j^* * \psi_j
\end {equation*}
holds in $L^p(G)$.
\end {proposition}
\begin {proof}
Since $ \displaystyle{\Sigma_n g = g * \underline{\delta_{2^{n+1}}^1} \phi - g * \underline{\delta_{2^{-n}}^1} \phi} $, Young's inequality implies that this sequence of operators is norm bounded. It is sufficient to show the convergence of the decomposition on the dense subspace $\mathcal{S}(G)$. We saw previously in Lemma~\ref{lemme_limite_at_infty} that $ g * \underline{\delta_{2^{n+1}}^1} \phi \stackrel[n \to +\infty]{}{\longrightarrow} c_{\phi} g $. Moreover, for $n \in \mathbb{N}$, one has :
\begin {equation*}
g * \underline{\delta_{2^{-n}}^1} \phi (x) = \frac{1}{2^{nQ}} \int_G g(y) \phi \bigg( \frac{1}{2^n} \odot (y^{-1} \cdot x) \bigg) dy = \frac{1}{2^{nQ}} \underline{\delta_{2^n}^1} g * \phi \bigg( \frac{1}{2^n} \odot x \bigg) \; \text{,}
\end {equation*}
and so :
\begin {equation*}
\| g * \underline{\delta_{2^{-n}}^1} \phi \|_{L^p(G)} = \frac{1}{2^{nQ}} \Bigg( \int_G \bigg| \underline{\delta_{2^n}^1} g * \phi \bigg( \frac{1}{2^n} \odot x \bigg) \bigg|^p dx \Bigg)^{\frac{1}{p}} = \frac{1}{2^{nQ\big( 1-\frac{1}{p} \big)}} \| \underline{\delta_{2^n}^1} g * \phi \|_{L^p(G)} \; \text{.}
\end {equation*}
Again, $ \underline{\delta_{2^n}^1} g * \phi \stackrel[n \to +\infty]{}{\longrightarrow} c_g \phi $, and in particular : $ \displaystyle{\frac{1}{2^{nQ\big( 1-\frac{1}{p} \big)}} \| \underline{\delta_{2^n}^1} g * \phi \|_{L^p(G)}} \stackrel[n \to +\infty]{}{\longrightarrow} 0 $. \\
Hence $ \Sigma_n g \stackrel[n \to +\infty]{}{\longrightarrow} c_{\phi} g $ and the case $p=2$ determines that $c_{\phi}=1$.
\end {proof}

\section {Characterization of Besov spaces by the discrete wavelets} \label{besov_spaces}
\noindent We show in this section that the characterization of $\mathring{B}^s_{p,q}(G)$ by a Littlewood-Paley theory can be discretized by sampling the convolution products $ f * \psi_j^* $ over a discrete set $ \Gamma \subset G $. This is equivalent to the study of the analysis operator $A_{\psi}$ : $ \mathcal{S}'(G) / \mathcal{P} \ni f \longmapsto A_{\psi}f = \big\{ \langle f,\psi_{j,\gamma} \rangle \big\}_{\begin {subarray}{l} j \in \mathbb{Z} \\ \gamma \in \Gamma \end {subarray}} $ associated to a discrete wavelet system $ \big\{ \psi_{j,\gamma} \big\}_{\begin {subarray}{l} j \in \mathbb{Z} \\ \gamma \in \Gamma \end {subarray}} $ defined by :
\begin {equation}
\forall x \in G, \; \psi_{j,\gamma}(x) = \underline{\tau_{\gamma}} \, \underline{\delta_{2^j}^1} \psi(x) = 2^{jQ} \psi(\gamma^{-1} \cdot 2^j \odot x) \; \text{,} \label{L^1_normalization}
\end {equation}
where $ \psi \in \mathcal{V}(G) $ is chosen as in Lemma~\ref{lemme_fondamental} and $\psi^*=\psi$. The main goal of this section is the proof of the equivalence between Besov norms and some associated discrete norm in Theorem~\ref{theoreme_echantillonnage}. \\
\newline
We cannot mention wavelets without some history. The development of wavelets started with \mbox{Alfred Haar} in 1909. Notable contributions can be attributed to George Zweig's discovery of the continuous wavelet transform (1975), Pierre Goupillaud, Alex Grossmann and Jean Morlet's formulation of the CWT (1982), Jan-Olov Str\"omberg's work on discrete wavelets (1983), Ingrid Daubechies' orthogonal wavelets with compact support (1988), Yves Meyer and St\'ephane Mallat's MRA framework (1989), and many more since. In 1986, Yves Meyer built a wavelet basis that suits for the simultaneous characterization of any Sobolev space $\mathring{H}^s$. Wavelet bases thus provide an elegant and unifying response to the problem of exhibiting unconditional bases for a wide range of classic functional spaces. Some useful references are \cite{Meyer1} and \cite{Meyer2}. \\
The very first wavelet bases on stratified Lie groups were obtained in 1989 by Pierre-Gilles Lemari\'e~\cite{Lemarie} from a spline interpolation theory. Already the sub-Laplacian $\Delta_G$ played a major role. Later in 2006 came out continuous and discrete wavelet systems by Daryl Geller and Azita Mayeli~\cite{GM} using the spectral theory of the sub-Laplacian. So it is rather natural to expect its involvement in a characterization of homogeneous Besov spaces. It was achieved in 2012 by Hartmut F\"uhr and Azita Mayeli~\cite{Fuhr}. In this section, we mainly rely on their analysis. \\
\newline \noindent In order to discretize norms over elementary tiles of $G$ and perform a multi resolution analysis (MRA), we shall introduce the notion of regular sampling sets.
\subsection {Regular sampling sets}
\begin {definition} \label{regular_sets}
Let $G$ be a Lie group. A subset $ \Gamma \subset G $ is called a regular sampling set if there is a relatively compact Borel neighborhood $ \mathcal{W} \subset G $ of the identity element of $G$ satisfying $ \displaystyle{\bigcup_{\gamma \in \Gamma} \gamma \mathcal{W} = G} $ (up to a set of measure zero), and for all $ \alpha, \gamma \in \Gamma, \; \alpha \neq \gamma, \; \mu_G( \alpha \mathcal{W} \cap \gamma \mathcal{W} ) = 0 $. Such a set $\mathcal{W}$ is called a $\Gamma$-tile.
\end {definition}
\begin {definition}
For $ \mathcal{U} \subset G $, a regular sampling set $\Gamma$ is said $\mathcal{U}$-dense if there exists a $\Gamma$-tile $ \mathcal{W} \subset \mathcal{U} $.
\end {definition}
\noindent For instance, lattices in $\mathbb{R}^d$ are a special class of regular sampling sets which are also cocompact discrete subgroups. However, some Lie groups fail to admit lattices, whereas by contrast there are always sufficiently dense regular sampling sets as indicated in the following lemma.
\begin {lemme}
Given a stratified Lie group $G$ with a Hausdorff geometrical realization, for any neighborhood $\mathcal{U}$ of the identity, there exists a $\mathcal{U}$-dense regular sampling set.
\end {lemme}
\begin {proof}
There exists $ \Gamma \subset G $ and a relatively compact set $\mathcal{W}$ with non-empty interior, such that $ \displaystyle{\bigcup_{\gamma \in \Gamma} \gamma \mathcal{W}} $ covers $G$ (possibly up to a set of measure zero). Then $ \mathcal{V} = \mathcal{W} x_0^{-1} $ is a $\Gamma$-tile for some point $x_0$ within $\mathcal{W}$. Finally, by choosing $\beta > 0$ small enough, we ensure that $\beta \mathcal{V} \subset \mathcal{U}$ and then $\beta \mathcal{V}$ is a $\beta \Gamma$-tile.
\end {proof}
\begin {definition}
Let $\Gamma \subset G$. An automorphism $\varrho \in Aut(G)$ is said $\Gamma$-acceptable if it leaves $\Gamma$ globally invariant, that is $ \varrho \Gamma \subseteq \Gamma $.
\end {definition}
\noindent We shall of course choose $\Gamma$ such that the dyadic dilations $\{\delta_{2^j}\}_{j \in \mathbb{Z}}$ and for any $ \gamma \in \Gamma $, the translations $\tau_{\gamma}$ are acceptable automorphisms. This definition ensures the compatibility between the group law, the nonisotropic dilations, the sampling set and the iterated wavelet system.
\subsection {Discretization of Besov norms}
\begin {definition}
Fix a discrete set $ \Gamma \subset G $. For any family $ \big\{ c_{j \gamma} \big\}_{\begin {subarray}{l} j \in \mathbb{Z} \\ \gamma \in \Gamma \end {subarray}} $ of complex numbers, we define :
\begin {equation}
\bigg\| \big\{ c_{j \gamma} \big\}_{\begin {subarray}{l} j \in \mathbb{Z} \\ \gamma \in \Gamma \end {subarray}} \bigg\|_{\mathring{b}^s_{p,q}} = \Bigg( \sum_{j \in \mathbb{Z}} \bigg( \sum_{\gamma \in \Gamma} \Big( 2^{j \big( s-\frac{Q}{p} \big)} |c_{j \gamma}| \Big)^p \bigg)^{\frac{q}{p}} \Bigg)^{\frac{1}{q}} \; \text{.}
\end {equation}
We introduce the space of coefficients $\mathring{b}^s_{p,q}(\Gamma)$ defined as :
\begin {equation}
\mathring{b}^s_{p,q}(\Gamma) = \Bigg\{ \big\{ c_{j \gamma} \big\}_{\begin {subarray}{l} j \in \mathbb{Z} \\ \gamma \in \Gamma \end {subarray}} \; \Bigg| \; \bigg\| \big\{ c_{j \gamma} \big\}_{\begin {subarray}{l} j \in \mathbb{Z} \\ \gamma \in \Gamma \end {subarray}} \bigg\|_{\mathring{b}^s_{p,q}} < +\infty \Bigg\} \; \text{,}
\end {equation}
which will be sometimes denoted by $\mathring{b}^s_{p,q}$ when there is no possible misunderstanding on $\Gamma$.
\end {definition}
\noindent The next theorem shows that Besov norms can be expressed in terms of discrete coefficients. The constants appearing in the norm equivalence depend on the functional spaces, but it does not matter if we use the same sampling set $\Gamma$ simultaneously for all spaces.
\begin {theoreme}[FM~\cite{Fuhr} - Theorem 5.4] \label{theoreme_echantillonnage}
There exists a neighborhood $\mathcal{U}$ of the identity such that, for any $\mathcal{U}$-dense regular sampling set $\Gamma$, one has :
\begin {equation*}
\forall u \in \mathcal{S}'(G) / \mathcal{P}, \; \forall 1 \leqslant p,q \leqslant +\infty, \; u \in \mathring{B}^s_{p,q}(G) \Longleftrightarrow \big\{ \langle u, 2^{-jQ} \psi_{j,\gamma} \rangle \big\}_{\begin {subarray}{l} j \in \mathbb{Z} \\ \gamma \in \Gamma \end {subarray}} \in \mathring{b}^s_{p,q}(\Gamma) \; \text{,}
\end {equation*}
where $\psi_{j,\gamma}$ is defined as in (\ref{L^1_normalization}). In addition, in $\mathring{B}^s_{p,q}(G)$ we have the norm equivalence :
\begin {equation} \label{equivalence_norm_besov}
\|u\|_{\mathring{B}^s_{p,q}(G)} \sim \Bigg( \sum_{j \in \mathbb{Z}} \bigg( \sum_{\gamma \in \Gamma} \Big( 2^{j \big( s-\frac{Q}{p} \big)} | \langle u, 2^{-jQ} \psi_{j,\gamma} \rangle | \Big)^p \bigg)^{\frac{q}{p}} \Bigg)^{\frac{1}{q}} \; \text{,}
\end {equation}
with some constants depending only on $p,q,s$ and the regular sampling set $\Gamma$.
\end {theoreme}
\noindent The accuracy of the estimate improves as the tightness of the sampling set increases.
\subsection {Banach wavelet frames for Besov spaces}
\noindent \newline In Hilbert spaces, a norm equivalence such as (\ref{equivalence_norm_besov}) is sufficient to conclude that the wavelet basis is a frame, allowing the reconstruction of $u$ from the discrete coefficients data. For Banach spaces, we need an extended definition of frames (refer to~\cite{Grochenig}), in order to show the invertibility of the associated frame operator. \\
After stating several technical lemmas, we will prove that any linear wavelet combination converges unconditionally in any $L^p(G)$ in Theorem~\ref{theoreme_inconditionnalite_lebesgue}, and also in $\mathring{B}^s_{p,q}(G)$ in Theorem~\ref{theoreme_inconditionnalite_besov} whenever its coefficients lie in $\mathring{b}^s_{p,q}$. Then for a suitable choice of sufficiently dense regular sampling sets $\Gamma$, the wavelet system $ \big\{ 2^{-jQ} \psi_{j,\gamma} \big\}_{\begin {subarray}{l} j \in \mathbb{Z} \\ \gamma \in \Gamma \end {subarray}} $ is a Banach frame for both $L^p(G)$ and $\mathring{B}^s_{p,q}(G)$.
\begin {definition}
A basis $(f_n)$ is an unconditional basis in a Banach space if any convergent series $ \displaystyle{\sum_n a_n f_n} $ converges unconditionally, that is to say the series $ \displaystyle{\sum_n a_{\sigma(n)} f_{\sigma(n)}} $ converges to the same limit for any \mbox{permutation $\sigma$} of the indexes.
\end {definition}
\noindent Recall that the sampled convolution products can be interpreted as scalar products too, that is to assert :
\begin {equation*}
f * \psi_j^* (2^{-j} \odot \gamma) = \langle f,\psi_{j,\gamma} \rangle \; \text{,}
\end {equation*}
where $\psi_{j,\gamma}$ denotes the wavelet of resolution $ \displaystyle{\frac{1}{2^j}} $ and at position $ \displaystyle{\frac{1}{2^j} \odot \gamma} $ - see (\ref{L^1_normalization}). The wavelet system is now used for synthesis.
\subsubsection {Unconditionality in $L^p(G)$ with $1<p<+\infty$}
\noindent \newline To check the unconditionality of wavelet decompositions in both Lebesgue and homogeneous Besov spaces, we need the next preliminary estimate.
\begin {lemme} \label{lemme_bound_column}
Let $ \eta,j \in \mathbb{Z} $ with $ \eta \leqslant j $ and $ n \geqslant Q+1 $. Let $\Gamma \subset G$ be a dense regular sampling set, in the sense of Definition~\ref{regular_sets}. Then :
\begin {equation*}
\forall x \in G, \; \sum_{\gamma \in \Gamma} \frac{2^{-jQ}}{\big( 1 + 2^\eta | 2^{-j} \odot \gamma^{-1} \cdot x |_G \big)^n} \leqslant c \; 2^{-\eta Q} \; \text{,}
\end {equation*}
where the constant $c$ depends only on $n$ and $\Gamma$.
\end {lemme}
\begin {proof}
By assumption, there exists a relatively compact open $\mathcal{W}$ such that $ \mu_G(\gamma \mathcal{W} \cap \gamma' \mathcal{W}) = 0 $ when $ \gamma \neq \gamma' $ in $\Gamma$. Then, in view of the left-invariance of the Haar measure $\mu_G$ in (\ref{left_invariance}), we have :
\begin {equation*}
\sum_{\gamma \in \Gamma} \frac{2^{-jQ}}{\big( 1 + 2^\eta | 2^{-j} \odot \gamma^{-1} \cdot x |_G \big)^n} \leqslant \sum_{\gamma \in \Gamma} \frac{1}{|\mathcal{W}|} \int_{2^{-j} \odot (\gamma \mathcal{W})} \frac{1}{\big( 1 + 2^\eta | 2^{-j} \odot \gamma^{-1} \cdot x |_G \big)^n} dy \; \text{.}
\end {equation*}
For any $y \in 2^{-j} \odot (\gamma \mathcal{W})$, the triangle inequality gives :
\begin {eqnarray*}
1 + 2^\eta | y^{-1} \cdot x |_G & \leqslant & 1 + 2^\eta c' \big( | y^{-1} \cdot 2^{-j} \odot \gamma |_G + | 2^{-j} \odot \gamma^{-1} \cdot x |_G \big) \\
& \leqslant & 1 + 2^\eta c' \big( 2^{-j} \text{ diam}(\mathcal{W}) + | 2^{-j} \odot \gamma^{-1} \cdot x |_G \big) \\
& \leqslant & c'' \big( 1 + 2^\eta | 2^{-j} \odot \gamma^{-1} \cdot x |_G \big)
\end {eqnarray*}
by using the fact that $\eta \leqslant j$. As a result :
\begin {eqnarray*}
\sum_{\gamma \in \Gamma} \frac{1}{|\mathcal{W}|} \int_{2^{-j} \odot (\gamma \mathcal{W})} \frac{1}{\big( 1 + 2^\eta | 2^{-j} \odot \gamma^{-1} \cdot x |_G \big)^n} dy & \leqslant & \bigg( \frac{1}{c''} \bigg)^n \sum_{\gamma \in \Gamma} \frac{1}{|\mathcal{W}|} \int_{2^{-j} \odot (\gamma \mathcal{W})} \frac{1}{\big( 1 + 2^\eta | y^{-1} \cdot x |_G \big)^n} dy \\
& = & 2^{-\eta Q} \bigg( \frac{1}{c''} \bigg)^n \frac{1}{|\mathcal{W}|} \int_G \frac{1}{\big( 1 + |y|_G \big)^n} dy \; \text{,}
\end {eqnarray*}
since the $\gamma \mathcal{W}$'s are pairwise disjoint. For $n \geqslant Q+1$, the latter integral is finite. Hence the lemma.
\end {proof}
\begin {theoreme} \label{theoreme_inconditionnalite_lebesgue}
Let $ 1 \leqslant p \leqslant +\infty $. Let $ \eta,j \in \mathbb{Z} $ be fixed with $ \eta \leqslant j $. Suppose that $ \Gamma \subset G $ is a regular sampling set. For all $ \gamma \in \Gamma $, we consider functions $f_{j,\gamma}$ on $G$ satisfying the following decay condition :
\begin {equation*}
\forall x \in G, \; \forall \eta,j \in \mathbb{Z}, \; \forall \gamma \in \Gamma, \; |f_{j,\gamma}(x)| \leqslant \frac{c_1}{\big( 1 + 2^\eta | 2^{-j} \odot \gamma^{-1} \cdot x |_G \big)^{Q+1}} \; \text{,}
\end {equation*}
with some constant $c_1>0$. We define $ \displaystyle{f_j = \sum_{\gamma \in \Gamma} c_{j \gamma} f_{j,\gamma}} $ with $ \{c_{j \gamma}\}_{\gamma} \in \ell^p(\Gamma) $. Then the series converges unconditionally in $L^p(G)$ with :
\begin {equation}
\|f_j\|_{L^p(G)} \leqslant c_2 \; 2^{(j-\eta)Q} 2^{- \frac{jQ}{p}} \bigg( \sum_{\gamma \in \Gamma} |c_{j \gamma}|^p \bigg)^{\frac{1}{p}} \; \text{,}
\end {equation}
for some constant $c_2$ independent of $j,\eta,\gamma$ and the sequence of coefficients $\{c_{j \gamma}\}_{\gamma}$.
\end {theoreme}
\begin {proof}
There exists a $\Gamma$-tile $\mathcal{W}$ such that $ \displaystyle{G = \bigsqcup_{\alpha \in \Gamma} 2^{-j} \odot (\alpha \mathcal{W})} $ . Then :
\begin {eqnarray*}
\|f_j\|_{L^p(G)}^p & = & \sum_{\alpha \in \Gamma} \int_{2^{-j} \odot (\alpha \mathcal{W})} \bigg| \sum_{\gamma \in \Gamma} c_{j \gamma} f_{j,\gamma}(x) \bigg|^p dx \\
& \leqslant & c_1^p \sum_{\alpha \in \Gamma} \int_{2^{-j} \odot (\alpha \mathcal{W})} \Bigg| \sum_{\gamma \in \Gamma} |c_{j \gamma}| \frac{1}{\big( 1 + 2^\eta | 2^{-j} \odot \gamma^{-1} \cdot x |_G \big)^{Q+1}} \Bigg|^p dx \; \text{.}
\end {eqnarray*}
On each integration domain $2^{-j} \odot (\alpha \mathcal{W})$, the triangle inequality yields the estimate :
\begin {equation*}
1 + 2^\eta | 2^{-j} \odot (\gamma^{-1} \cdot \alpha) |_G \leqslant c'' \big( 1 + 2^\eta | 2^{-j} \odot \gamma^{-1} \cdot x |_G \big) \; \text{.}
\end {equation*}
Then the integrand can be majorized by the constant :
\begin {equation*} \Bigg| \sum_{\gamma \in \Gamma} |c_{j \gamma}| \frac{c''}{\big( 1 + 2^\eta | 2^{-j} \odot (\gamma^{-1} \cdot \alpha) |_G \big)^{Q+1}} \Bigg|^p \; \text{.}
\end {equation*}
Hence :
\begin {eqnarray*}
\|f_j\|_{L^p(G)}^p & \leqslant & (c_1 c'')^p |\mathcal{W}| \sum_{\alpha \in \Gamma} 2^{-jQ} \Bigg( \sum_{\gamma \in \Gamma} |c_{j \gamma}| \frac{1}{\big( 1 + 2^\eta | 2^{-j} \odot (\gamma^{-1} \cdot \alpha) |_G \big)^{Q+1}} \Bigg)^p \\
& = & (c_1 c'')^p |\mathcal{W}| \sum_{\alpha \in \Gamma} 2^{-jQ} \bigg( \sum_{\gamma \in \Gamma} |c_{j \gamma}| a_{\alpha \gamma} \bigg)^p \; \text{,}
\end {eqnarray*}
with $ \displaystyle{a_{\alpha \gamma} = \frac{1}{\big( 1 + 2^\eta | 2^{-j} \odot (\gamma^{-1} \cdot \alpha) |_G \big)^{Q+1}}}  \; \cdot $ Lemma~\ref{lemme_bound_column} now ensures that the Schur lemma's conditions are fulfilled for the coefficients $\{a_{\alpha \gamma}\}$ with $ \displaystyle{\max \bigg( \sup_\alpha \sum_{\gamma \in \Gamma} |a_{\alpha \gamma}|, \sup_\gamma \sum_{\alpha \in \Gamma} |a_{\alpha \gamma}| \bigg) \leqslant 2^{(j-\eta)Q}} $. Therefore, there exists a constant $c_2$ such that :
\begin {equation*}
\|f_j\|_{L^p(G)} \leqslant 2^{-\frac{jQ}{p}} c_1 c'' |\mathcal{W}|^{\frac{1}{p}} \Bigg( \sum_{\alpha \in \Gamma} \bigg( \sum_{\gamma \in \Gamma} |c_{j \gamma}| a_{\alpha \gamma} \bigg)^p \Bigg)^\frac{1}{p} \leqslant 2^{-\frac{jQ}{p}} 2^{(j-\eta)Q} c_2 \bigg( \sum_{\gamma \in \Gamma} |c_{j \gamma}|^p \bigg)^{\frac{1}{p}} \; \text{.}
\end {equation*}
\end {proof}
\subsubsection {Unconditionality in $\mathring{B}^s_{p,q}(G)$ with $1 \leqslant p,q < +\infty$}
\begin {lemme}
With the above notations, there exists a constant $c>0$ such that $ \forall j,\ell \in \mathbb{Z}, \; \forall \gamma \in \Gamma, \; \forall x \in G $, the following estimate holds :
\begin {equation} \label{bound_convolution}
| \psi_{j,\gamma} * \psi_{\ell}^* (x) | \leqslant \left\{
\begin {array}{cl}
\displaystyle{\frac{c \; 2^{jQ}}{\big( 1 + 2^j | 2^{-j} \odot \gamma^{-1} \cdot x |_G \big)^{Q+1}}} & \text{if } |\ell-j| \leqslant 1 \\
0 & \text{otherwise}
\end {array} \right. \; \text{.}
\end {equation}
\end {lemme}
\begin {proof}
Let us compute :
\begin {eqnarray*}
\psi_{j,\gamma} * \psi_{\ell}^* (x) & = & \int_G 2^{jQ} \psi(\gamma^{-1} \cdot 2^j \odot y) 2^{\ell Q} \overline{\psi\big(2^{\ell} \odot (x^{-1} \cdot y)\big)} dy \\
& = & \int_G 2^{jQ} \psi(y) \overline{\psi_{\ell-j}\big((\gamma^{-1} \cdot 2^j \odot x)^{-1} \cdot y)\big)} dy \\
& = & 2^{jQ} ( \psi * \psi_{\ell-j}^* )(\gamma^{-1} \cdot 2^j \odot x) \; \text{.}
\end {eqnarray*}
By (\ref{relation_rangs}), this convolution product vanishes when $ |j-\ell| > 1 $. \\
In the other case i.e. $ \ell-j \in \{-1,0,+1\} $, the convolution products $ \psi * \psi_{\ell-j} $ are functions in the Schwartz class and henceforth, $ | \psi_{j,\gamma} * \psi_{\ell}^* (x) | \leqslant \displaystyle{\frac{c \; 2^{jQ}}{\big( 1 + 2^j | 2^{-j} \odot \gamma^{-1} \cdot x |_G \big)^{Q+1}}} $ for some constant $c$.
\end {proof}
\noindent The next result is crucial for the upcoming Section~\ref{profils_stratified_groups}.
\begin {theoreme} \label{theoreme_inconditionnalite_besov}
Let $ 1 \leqslant p,q < +\infty $. If $\Gamma$ is the same regular sampling set as in Theorem~\ref{theoreme_inconditionnalite_lebesgue}, then for any sequence of coefficients $ \big\{ c_{j \gamma} \big\}_{\begin {subarray}{l} j \in \mathbb{Z} \\ \gamma \in \Gamma \end {subarray}} \in \mathring{b}^s_{p,q}(\Gamma) $, the sum $ \displaystyle{f = \sum_{j,\gamma} 2^{-jQ} c_{j \gamma} \psi_{j,\gamma}} $ converges unconditionally in Besov norm with :
\begin {equation}
\|f\|_{\mathring{B}^s_{p,q}(G)} \leqslant c \; \Bigg( \sum_{j \in \mathbb{Z}} \bigg( \sum_{\gamma \in \Gamma} \Big( 2^{j \big( s-\frac{Q}{p} \big)} |c_{j \gamma}| \Big)^p \bigg)^{\frac{q}{p}} \Bigg)^{\frac{1}{q}} \; \text{,}
\end {equation}
for some constant $c$ independent of $ \big\{ c_{j \gamma} \big\}_{\begin {subarray}{l} j \in \mathbb{Z} \\ \gamma \in \Gamma \end {subarray}} $.
\end {theoreme}
\begin {proof}
It is sufficient to prove the norm estimate for finite sequences of coefficients $ \big\{ c_{j \gamma} \big\}_{\begin {subarray}{l} j \in \mathbb{Z} \\ \gamma \in \Gamma \end {subarray}} $. The full statement follows from the completeness of $\mathring{B}^s_{p,q}(G)$ and the property that the Kronecker symbols $\delta$ form an unconditional basis of $\mathring{b}^s_{p,q}$ (this is why $p,q < +\infty $ is required). By definition :
\begin {eqnarray*}
\|f\|_{\mathring{B}^s_{p,q}(G)} & = & \Big\| \big\{ 2^{\ell s}\| f * \psi_{\ell}^* \|_{L^p(G)} \big\}_{\ell \in \mathbb{Z}} \Big\|_{\ell^q(\mathbb{Z})} = \Bigg\| \Bigg\{ 2^{\ell s} \bigg\| \sum_{j,\gamma} 2^{-jQ} c_{j \gamma} \psi_{j,\gamma} * \psi_{\ell}^* \bigg\|_{L^p(G)} \Bigg\}_{\ell \in \mathbb{Z}} \Bigg\|_{\ell^q(\mathbb{Z})} \\
& = & \Bigg\| \Bigg\{ 2^{\ell s} \bigg\| \sum_{\gamma \in \Gamma} \sum_{j= \ell -1}^{\ell+1} 2^{-jQ} c_{j \gamma} \psi_{j,\gamma} * \psi_{\ell}^* \bigg\|_{L^p(G)} \Bigg\}_{\ell \in \mathbb{Z}} \Bigg\|_{\ell^q(\mathbb{Z})} \\
& \leqslant & \Bigg\| \Bigg\{ 2^{\ell s} \sum_{j= \ell -1}^{\ell+1} \bigg\| \sum_{\gamma \in \Gamma} 2^{-jQ} c_{j \gamma} \psi_{j,\gamma} * \psi_{\ell}^* \bigg\|_{L^p(G)} \Bigg\}_{\ell \in \mathbb{Z}} \Bigg\|_{\ell^q(\mathbb{Z})}
\end {eqnarray*}
by using (\ref{relation_rangs}). For instance, let us consider the middle term when $j=\ell$. Applying successively (\ref{bound_convolution}) and Theorem~\ref{theoreme_inconditionnalite_lebesgue} where $j-\eta=0$, we obtain that :
\begin {equation*}
2^{\ell s} \bigg\| \sum_{\gamma \in \Gamma} 2^{-\ell Q} c_{\ell \gamma} \psi_{\ell,\gamma} * \psi_{\ell}^* \bigg\|_{L^p(G)} \leqslant 2^{\ell \big( s - \frac{Q}{p} \big)} c_2 \bigg( \sum_{\gamma \in \Gamma} |c_{\ell \gamma}|^p \bigg)^{\frac{1}{p}} \; \text{.}
\end {equation*}
Therefore, we get :
\begin {equation*}
\Bigg\| \Bigg\{ 2^{\ell s} \bigg\| \sum_{\gamma \in \Gamma} 2^{-\ell Q} c_{\ell \gamma} \psi_{\ell,\gamma} * \psi_{\ell}^* \bigg\|_{L^p(G)} \Bigg\}_{\ell \in \mathbb{Z}} \Bigg\|_{\ell^q(\mathbb{Z})} \leqslant c_2 \; \Bigg( \sum_{\ell \in \mathbb{Z}} \bigg( \sum_{\gamma \in \Gamma} \Big( 2^{\ell \big( s-\frac{Q}{p} \big)} |c_{\ell \gamma}| \Big)^p \bigg)^{\frac{q}{p}} \Bigg)^{\frac{1}{q}} = c_2 \; \bigg\| \big\{ c_{j \gamma} \big\}_{\begin {subarray}{l} j \in \mathbb{Z} \\ \gamma \in \Gamma \end {subarray}} \bigg\|_{\mathring{b}^s_{p,q}} \; \text{.}
\end {equation*}
Of course, by applying similar calculations when $j=\ell \pm 1$, we will find that $ \|f\|_{\mathring{B}^s_{p,q}(G)} \leqslant c \; \bigg\| \big\{ c_{j \gamma} \big\}_{\begin {subarray}{l} j \in \mathbb{Z} \\ \gamma \in \Gamma \end {subarray}} \bigg\|_{\mathring{b}^s_{p,q}} $.
\end {proof}
\noindent Generally, for sufficiently dense regular sampling sets $\Gamma$, the invertibility of the synthesis operator \mbox{$ S_{\psi,\Gamma} $ : $ f \longmapsto S_{\psi,\Gamma} f = \displaystyle{\sum_{j,\gamma} 2^{-jQ} \langle f, \psi_{j,\gamma} \rangle \psi_{j,\gamma}} $} ensures that the discrete wavelet basis obtained by dilations by powers of 2 and translations provide an universal Banach frame simultaneously for all Besov spaces $\mathring{B}^s_{p,q}(G)$ with $s \in \mathbb{R}$ and $ 1 \leqslant p,q < +\infty $.

\section {Construction of the profiles on stratified Lie groups} \label{profils_stratified_groups}
\noindent At this stage, we have enough tools to prove Theorem~\ref{theoreme_profils} by adapting the method in~\cite{BCK}. To sum up, we have constructed LP-admissible wavelets $ \psi \in \mathcal{V}(G) $, and then confirmed the existence of unconditional bases (spanned as an iterated functions system by $\psi$) \mbox{for both $L^p(G)$ and $\mathring{B}^s_{p,q}(G)$ in Theorems~\ref{theoreme_inconditionnalite_lebesgue} and \ref{theoreme_inconditionnalite_besov}.} \\
\newline Let $\Gamma \subset G$ be a regular sampling set such that $\{\tau_{\gamma}\}_{\gamma \in \Gamma}$ and $\{\delta_{2^j}\}_{j \in \mathbb{Z}}$ are \mbox{$\Gamma$-acceptable automorphisms.} We introduce a new index $\lambda = (j,\gamma)$ for the wavelet basis, where the component $j=j_\lambda$ is a scale index (actually the dyadic exponent) while the other component $\gamma = \gamma_\lambda$ is a space index. For any element $f$ in $ \mathring{H}^s(G) = \mathring{B}^s_{2,2}(G) $, its wavelet decomposition can be written as :
\begin {equation*}
f = \sum_{\lambda \in \mathbb{Z} \times \Gamma} d_\lambda \psi_\lambda \; \text{,}
\end {equation*}
where $\psi_{\lambda}$ is now $L^p$-normalized according to :
\begin {equation}
\forall \lambda \in \mathbb{Z} \times \Gamma, \; \psi_{\lambda}(x) = 2^{\frac{j_\lambda Q}{p}} \psi(\gamma_{\lambda}^{-1} \cdot 2^{j_\lambda} \odot x) = \underline{\tau_{\gamma_{\lambda}}} \, \underline{\delta_{2^{j_\lambda}}^p} \psi(x) \; \text{.} \label{scaling}
\end {equation}
By Theorem~\ref{psi_independent}, the homogeneous Besov spaces are independent of the choice of $\psi \in \mathcal{V}(G)$. By picking $\hat{\psi}$ with support in the interval $ \displaystyle{\bigg[ \frac{1}{2},1 \bigg]} $, the $\psi_\lambda$'s shall form an orthonormal basis. From~(\ref{scaling}), one has readily :
\begin {eqnarray}
& & \left\{
\begin {array}{l}
\| \psi_{\lambda} \|_{L^p(G)} = \| \psi \|_{L^p(G)} \\
\| \psi_{\lambda} \|_{\mathring{H}^s(G)} = \| \psi \|_{\mathring{H}^s(G)}
\end {array} \right. \text{,} \label{norm_invariance_by_scaling} \\
\text{as well as} & & \|f\|_{\mathring{H}^s(G)} = \| \{d_\lambda\} \|_{\ell^2(\mathbb{Z} \times \Gamma)} \; \text{.}
\end {eqnarray}
The unconditionality of this basis implies the existence of a constant $D$ such that for any finite subset $ E \subset \mathbb{Z} \times \Gamma $, any coefficients $(c_\lambda)_{\lambda \in E}$ and $(d_\lambda)_{\lambda \in E}$ satisfying $ \forall \lambda, \; |c_\lambda| \leqslant |d_\lambda| $, one has :
\begin {equation}
\bigg\| \sum_{\lambda \in E} c_\lambda \psi_\lambda \bigg\|_{\mathring{H}^s(G)} \leqslant D \; \bigg\| \sum_{\lambda \in E} d_\lambda \psi_\lambda \bigg\|_{\mathring{H}^s(G)} \; \text{.} \label{propriete_inconditionnalite}
\end {equation}
\newline
For $M>0$, let us now consider the nonlinear projector $Q_M$ defined by :
\begin {equation}
\forall f \in \mathring{H}^s(G), \; Q_M(f) = \sum_{\lambda \in E_M(f)} d_\lambda \psi_\lambda \; \text{,} \label{nonlinear_projection}
\end {equation}
where $E_M(f)$ is the subset of $ \mathbb{Z} \times \Gamma $ of cardinality $M$ corresponding to the $M$ largest values of $|d_\lambda|$. Note that such a set always exists, but the nonlinear projection~(\ref{nonlinear_projection}) may not be unique when some $|d_{\lambda}|$ are equal. In which case, any realization of such set suits for $E_M(f)$. Additionally, it has been proven in~\cite{DeVore} or~\cite{DVJP} that :
\begin {eqnarray}
\lim_{M \to +\infty} \sup_{\|f\|_{\mathring{H}^s(G)} \leqslant 1} \| f-Q_M(f) \|_{L^p(G)} = 0 \; \text{,} \label{uniform_convergence_Q_M}
\end {eqnarray}
with $ \displaystyle{\frac{s}{Q} + \frac{1}{p} = \frac{1}{2}}  \; \cdot $ Let us emphasize that the uniform convergence~(\ref{uniform_convergence_Q_M}) of $Q_M(f)$ to $f$ in $L^p(G)$ is tied to the nonlinear nature of the operator $Q_M$. The nonlinear projection $Q_M(f)$, sometimes called the best $M$-term approximation of $f$, has been extensively studied - refer for instance to Ronald A. DeVore~\cite{DeVore} and the numerous references therein. \\
\newline
The proof of Theorem~\ref{theoreme_profils} is based on a diagonal subsequence, and it is structured in three main steps. We work of course under the theorem's assumptions and consider a sequence $(u_n)_{n>0}$ of bounded functions in $\mathring{H}^s(G)$. Then let us define :
\begin {equation*}
K = \sup_{n>0} \| u_n \|_{\mathring{H}^s(G)} < +\infty \; \text{.}
\end {equation*}
\subsection {Reordering of the wavelet decomposition} \label{reorder_wavelet_coefficients}
\noindent \newline From the wavelet decomposition $ \displaystyle{u_n = \sum_{\lambda \in \mathbb{Z} \times \Gamma} d_{\lambda,n} \psi_\lambda} $, the summands are reordered by decreasing moduli $|d_{\lambda,n}|_{\lambda \in \mathbb{Z} \times \Gamma}$ such that :
\begin {equation*}
u_n = \sum_{m>0} d_{m,n} \psi_{\lambda(m,n)} \; \text{.}
\end {equation*}
Using the nonlinear projector defined by~(\ref{nonlinear_projection}), one gets :
\begin {equation*}
u_n = \sum_{m=1}^M d_{m,n} \psi_{\lambda(m,n)} + \big( u_n - Q_M(u_n) \big) \; \text{,}
\end {equation*}
with, in light of~(\ref{uniform_convergence_Q_M}) :
\begin {equation}
\lim_{M \to +\infty} \sup_{n>0} \| u_n - Q_M(u_n) \|_{L^p(G)} = 0 \; \text{.} \label{uniform_convergence_M_approximation}
\end {equation}
Since wavelets are normalized in $ \mathring{H}^s(G) $, we know that : $ \displaystyle{\sup_{m \geqslant 1,n} |d_{m,n}| \leqslant DK} $,
where $D$ is the constant of~(\ref{propriete_inconditionnalite}). Up to a possible diagonal extraction in $n$ of a subsequence, we can assume that, for $m \geqslant 1$, the sequence $(d_{m,n})_{n>0}$ converges to a finite limit depending only on $m$ :
\begin {equation*}
d_m = \lim_{n \to +\infty} d_{m,n} \; \text{.}
\end {equation*}
We can then write :
\begin {equation}
u_n = \sum_{m=1}^M d_m \psi_{\lambda(m,n)} + \sum_{m=1}^M (d_{m,n} - d_m) \psi_{\lambda(m,n)} + \big( u_n - Q_M(u_n) \big) \; \text{.} \label{asymptotic_expansion}
\end {equation}
\subsection {Extraction of the approximate profiles} \label{extraction_approximate_profiles}
\noindent \newline The exact profiles $\phi^{\ell}$ involved in~(\ref{profils}) are inferred as limits in $\mathring{H}^s(G)$ of some approximate profiles $\phi^{\ell,i}$ obtained by the following procedure :
\begin {enumerate}
\item Initialize $ \phi^{1,1} = d_1 \psi $, $ \lambda_1(n) = \lambda(1,n) $ and $ \varphi_1(n) = n $.
\item At step $i-1$, assume that we have obtained $\nu(i-1)$ functions ($\phi^{1,i-1}$, $\phi^{2,i-1}$, \ldots, $\phi^{\nu(i-1),i-1}$), scale-space indexes $(\lambda_1(n)$, $\lambda_2(n)$, \ldots, $\lambda_{\nu(i-1)}(n)$), as well as an increasing sequence of positive integers $\varphi_{i-1}(n)$ such that :
\begin {equation*}
\sum_{m=1}^{i-1} d_m \psi_{\lambda(m,\varphi_{i-1}(n))} = \sum_{\ell=1}^{\nu(i-1)} \phi_{\lambda_{\ell}(\varphi_{i-1}(n))}^{\ell,i-1} \; \text{,}
\end {equation*}
where $ \displaystyle{\phi_{\lambda_\ell}^{\ell,i-1} = \underline{\tau_{\gamma_{\lambda_\ell}}} \, \underline{\delta_{2^{j_{\lambda_\ell}}}^p} \phi^{\ell,i-1}} $ as in (\ref{scaling}). Superscripts in $\phi_{\lambda_\ell}^{\ell,i-1}$ are harmless summation indexes, but the subscript indicates a translated and dilated copy of $\phi^{\ell,i-1}$.
\item Add the $i$-th term $d_i \psi_{\lambda(i,\varphi_{i-1}(n))}$ either to build a new function, either to modify slightly one of the previous functions according to the next dichotomy :
\begin {itemize}
\item[\textbullet] Case 1 : Assume that we can extract $\varphi_i(n)$ from $\varphi_{i-1}(n)$ such that for any $ \ell \in \llbracket 1,\nu(i-1) \rrbracket  $, at least one of the two conditions below is satisfied :
\begin {eqnarray*}
& & \lim_{n \to +\infty} \big| j(\lambda_{\ell}(\varphi_i(n))) - j(\lambda(i,\varphi_i(n))) \big| = + \infty \; \text{,} \\
\text{or} & & \lim_{n \to +\infty} \bigg| \frac{2^{j(\lambda(i,\varphi_i(n)))}}{2^{j(\lambda_{\ell}(\varphi_i(n)))}} \odot \gamma(\lambda_{\ell}(\varphi_i(n)))^{-1} \cdot \gamma(\lambda(i,\varphi_i(n))) \bigg|_G = +\infty \; \text{.}
\end {eqnarray*}
In that case, with $\nu(i)=\nu(i-1)+1$, we add a new function $\phi^{\nu(i),i}$ such that :
\begin {eqnarray*}
\phi^{\nu(i),i} = d_i \psi, \; \lambda_{\nu(i)}(n) = \lambda(i,n) \; \text{,}
\end {eqnarray*}
and we keep every previous approximate profiles, by setting : $ \forall \ell \in \llbracket 1,\nu(i-1) \rrbracket  $, $ \phi^{\ell,i} = \phi^{\ell,i-1} $.
\item[\textbullet] Case 2 : Suppose that for some subsequence $\varphi_i(n)$ of $\varphi_{i-1}(n)$ and some $\ell \in \llbracket 1,\nu(i-1) \rrbracket $, none of the two above conditions is satisfied. Then one can check that $ j(\lambda_{\ell}(\varphi_i(n))) - j(\lambda(i,\varphi_i(n))) $ and $ \displaystyle{\bigg| \frac{2^{j(\lambda(i,\varphi_i(n)))}}{2^{j(\lambda_{\ell}(\varphi_i(n)))}} \odot \gamma(\lambda_{\ell}(\varphi_i(n)))^{-1} \cdot \gamma(\lambda(i,\varphi_i(n))) \bigg|_G} $ only take a finite number of values as $n$ varies. Therefore, up to an additional subsequence extraction, we can infer the existence of finite values $\tilde{\jmath} \in \mathbb{Z}$ and $\tilde{\gamma} \in G$ such that $\forall n>0$ :
\begin {eqnarray*}
& & j(\lambda(i,\varphi_i(n))) - j(\lambda_{\ell}(\varphi_i(n))) = \tilde{\jmath} \; \text{,} \\
\text{and} & & \frac{2^{j(\lambda(i,\varphi_i(n)))}}{2^{j(\lambda_{\ell}(\varphi_i(n)))}} \odot \gamma(\lambda_{\ell}(\varphi_i(n)))^{-1} \cdot \gamma(\lambda(i,\varphi_i(n))) = \tilde{\gamma} \; \text{.}
\end {eqnarray*}
Set $\nu(i)=\nu(i-1)$. The function $\phi^{\ell,i-1}$ is now replaced by :
\begin {equation*}
\phi^{\ell,i}(x) = \phi^{\ell,i-1}(x) + 2^{\frac{\tilde{\jmath} Q}{p}} d_i \psi \big( \tilde{\gamma}^{-1} \cdot 2^{\tilde{\jmath}} \odot x \big) \; \text{,}
\end {equation*}
whereas the other profiles remain unchanged i.e. for all $ \llbracket 1,\nu(i-1) \rrbracket  \ni \ell' \neq \ell$, $ \phi^{\ell',i} = \phi^{\ell',i-1} $.
\end {itemize}
\end {enumerate}
This algorithm shows that, for any $M \geqslant 1$, there exists $ \nu(M) \leqslant M $ such that :
\begin {equation*}
\sum_{m=1}^M d_m \psi_{\lambda(m,n)} = \sum_{\ell=1}^{\nu(M)} \phi_{\lambda_{\ell}(n)}^{\ell,M} \; \text{.}
\end {equation*}
More explicitly, for every $ \ell \in \llbracket 1,\nu(M)\rrbracket  $, we have :
\begin {equation*}
\phi_{\lambda_{\ell}(n)}^{\ell,M} = \sum_{m \in E(\ell,M)} d_m \psi_{\lambda(m,n)} \; \text{,}
\end {equation*}
where the sets $E(\ell,M)$ form a disjoint partition of $ \displaystyle{\llbracket 1,M \rrbracket  = \bigsqcup_{\ell=1}^{\nu(M)} E(\ell,M)} $. \\
It is clear that $ E(\ell,M) \subseteq E(\ell,M+1) $, and the number of approximate profiles $\nu(M)$ increases by at most one unit when going from $M$ to $M+1$.
\subsection {End of the proof}
\noindent \newline To finish the proof of Theorem~\ref{theoreme_profils}, the exact profiles $\phi^{\ell}$ are obtained as limits in $\mathring{H}^s(G)$ of the approximate profiles $\phi^{\ell,M}$ as $M \to +\infty$, the same way as in~\cite{BCK} by means of the invariance by scaling~(\ref{norm_invariance_by_scaling}) and the unconditionality of the wavelet basis. \\
\newline
Let us now estimate the error terms in (\ref{asymptotic_expansion}). \\
For a given $ L \in \llbracket 1,M \rrbracket $, according to Subsections~\ref{reorder_wavelet_coefficients} and~\ref{extraction_approximate_profiles}, we can rewrite $u_n$ as :
\begin {equation}
u_n = \sum_{\ell=1}^L \phi_{\lambda_{\ell}(n)}^{\ell} + r_{n,L} \; \text{,}
\end {equation}
where the remainder $r_{n,L}$ can be split into :
\begin {equation*}
r_{n,L} = \underbrace{ \sum_{\ell=1}^L \Big( \phi_{\lambda_{\ell}(n)}^{\ell,M} - \phi_{\lambda_{\ell}(n)}^{\ell} \Big) + \sum_{\ell=1}^L \sum_{m \in E(\ell,M)} (d_{m,n} - d_m) \psi_{\lambda(m,n)} }_{r_1(n,L,M)} + \underbrace{\sum_{\ell=L+1}^{\nu(M)} \sum_{m \in E(\ell,M)} d_{m,n} \psi_{\lambda(m,n)} + \big( u_n - Q_M(u_n) \big) }_{r_2(n,L,M)} \; \text{.}
\end {equation*}
Observe that each of these summands depends on the chosen value of $M$, but their total sum $r_{n,L}$ is actually independent of $M$. \\
\newline
Under the norm invariance by scaling~(\ref{norm_invariance_by_scaling}) of the $L^p$-normalized $\{\psi_{\lambda}\}$ basis, we infer that
\begin {equation*}
\bigg\| \sum_{\ell=1}^L \Big( \phi_{\lambda_{\ell}(n)}^{\ell,M} - \phi_{\lambda_{\ell}(n)}^{\ell} \Big) \bigg\|_{\mathring{H}^s(G)} \leqslant \sum_{\ell=1}^L \bigg\| \Big( \phi_{\lambda_{\ell}(n)}^{\ell,M} - \phi_{\lambda_{\ell}(n)}^{\ell} \Big) \bigg\|_{\mathring{H}^s(G)} = \sum_{\ell=1}^L \big\| \phi^{\ell,M} - \phi^{\ell} \big\|_{\mathring{H}^s(G)} \; \text{.}
\end {equation*}
Since for all $\ell \geqslant 1$, $ \phi^{\ell,M} \stackrel[M \to +\infty]{}{\longrightarrow} \phi^{\ell} $ in $\mathring{H}^s(G)$, one deduces that for any fixed $L \geqslant 1$ :
\begin {equation*}
\varlimsup_{n\to +\infty} \bigg\| \sum_{\ell=1}^L \Big( \phi_{\lambda_{\ell}(n)}^{\ell,M} - \phi_{\lambda_{\ell}(n)}^{\ell} \Big) \bigg\|_{\mathring{H}^s(G)} 
 \stackrel[M \to +\infty]{}{\longrightarrow}
 0 \; \text{.}
\end {equation*}
Now combining (\ref{propriete_inconditionnalite}) and the norm invariance~(\ref{norm_invariance_by_scaling}), for all $M$ and $ 1 \leqslant L \leqslant \nu(M) $ fixed, one has :
\begin {eqnarray*}
\bigg\| \sum_{\ell=1}^L \sum_{m \in E(\ell,M)} (d_{m,n} - d_m) \psi_{\lambda(m,n)} \bigg\|_{\mathring{H}^s(G)} & \leqslant & D \; \bigg\| \sum_{m=1}^M (d_{m,n} - d_m) \psi_{\lambda(m,n)} \bigg\|_{\mathring{H}^s(G)} \\
& \leqslant & D \sum_{m=1}^M |d_{m,n} - d_m| \| \psi \|_{\mathring{H}^s(G)} \; \text{.}
\end {eqnarray*}
Consequently : $$ \forall L,M \geqslant 1, \; \bigg\| \sum_{\ell=1}^L \sum_{m \in E(\ell,M)} (d_{m,n} - d_m) \psi_{\lambda(m,n)} \bigg\|_{\mathring{H}^s(G)} \stackrel[n \to +\infty]{}{\longrightarrow} 0 \; \text{.} $$ So we get :
\begin {equation*}
\forall L \geqslant 1, \; \varlimsup_{n \to +\infty} \| r_1(n,L,M) \|_{\mathring{H}^s(G)} \stackrel[M \to +\infty]{}{\longrightarrow} 0 \; \text{,}
\end {equation*}
which in view of (\ref{critical_sobolev_injection}) ensures that the same holds for $ \| r_1(n,L,M) \|_{L^p(G)} $. \\
\newline
Moreover, the term $r_2(n,L,M)$ can be viewed as the partial sum $\displaystyle{\sum_{\ell \geqslant L+1} d_{m,n} \psi_{\lambda(m,n)}}$. So :
\begin {equation*}
\| r_2(n,L,M) \|_{L^p(G)} \leqslant D \; \bigg\| \sum_{\ell \geqslant L+1} d_{m,n} \psi_{\lambda(m,n)} \bigg\|_{L^p(G)} \; \text{,}
\end {equation*}
and by (\ref{uniform_convergence_M_approximation}), its convergence to $0$ when $L \to +\infty$ is assured. Since $M \geqslant L$, one obtains :
\begin{equation*}
\lim_{L \to +\infty} \sup_{n>0} \| r_2(n,L,M) \|_{L^p(G)} = 0 \; \text{.}
\end{equation*}
\newline
Lastly, the property :
\begin {eqnarray*}
& & \| u_n \|_{\mathring{H}^s(G)}^2 = \sum_{\ell=1}^L \| \phi^{\ell} \|_{\mathring{H}^s(G)}^2 + \| r_{n,L} \|_{\mathring{H}^s(G)}^2 + o(1) \text{ as } n \to +\infty \; \text{,} \\
& \text{with} & \lim_{L \to +\infty} \varlimsup_{n \to +\infty} \| r_{n,L} \|_{L^p(G)} = 0 \; \text{,}
\end {eqnarray*}
follows as a corollary of the wavelets' mutual orthogonality and their well-defined $L^2$-normalization in~(\ref{profils}). That eventually concludes Theorem~\ref{theoreme_profils}'s proof.
\begin {remarque}
As a final observation, note that due to multiple extractions of subsequences - and the use of the underlying axiom of choice (AC), it turns out that the profile decomposition may yet not be unique.
\end {remarque}

\noindent \newline

\end {document}